\documentclass[10pt, oneside, a4paper]{article}
\usepackage[titletoc]{appendix}
\usepackage{amsmath, amsthm, amssymb, amsfonts}
\usepackage{geometry}
\usepackage{graphicx}
\usepackage{natbib}
\usepackage{setspace}
\usepackage{algorithm}
\usepackage{algorithmicx}
\usepackage{algpseudocode}
\usepackage{longtable}
\usepackage{caption}
\usepackage{float}
\usepackage{enumerate}
\usepackage{url}

\newtheorem{theorem}{Theorem}[section]
\newtheorem{lemma}[theorem]{Lemma}
\newtheorem{proposition}[theorem]{Proposition}

\theoremstyle{definition}

\newtheorem{example}{Example}[section]

\newtheorem{remark}[theorem]{Remark}

\newcommand{\renum}[1]{}
\input{tcilatex}

\begin{document}

\title{A Modern Gauss-Markov Theorem? Really?\thanks{%
We would like to thank Abram Kagan for answering an inquiry and pointing us
to the reference \cite{KS}, and Peter Phillips for comments on the first
version of the paper. We are also grateful to Stephen Portnoy for sending us
his paper, and to Bruce Hansen for making the page proofs available to us
after we had sent him the first version of our paper. An abridged version of
the present paper has recently been accepted by Econometrica. Address
correspondence to Benedikt P\"{o}tscher, Department of Statistics,
University of Vienna, A-1090 Oskar-Morgenstern Platz 1. E-Mail:
benedikt.poetscher@univie.ac.at. }}
\author{Benedikt M. P\"{o}tscher and David Preinerstorfer \\
Department of Statistics, University of Vienna\\
SEPS-SEW, University of St.~Gallen}
\date{First version: February 2022\\
Second version: March 2022\\
Third version: May 2022\\
Fourth version: December 2022\\
Fifth version: October 2023\\
}
\maketitle

\begin{abstract}
We show that the theorems in \cite{Hansen_1.5} (the version accepted by
Econometrica), except for one, are not new as they coincide with classical
theorems like the good old Gauss-Markov or Aitken Theorem, respectively; the
exceptional theorem is incorrect. \cite{Hansen_2} corrects this theorem. As
a result, all theorems in the latter version coincide with the above
mentioned classical theorems. Furthermore, we also show that the theorems in 
\cite{Hansen_3} (the version published in Econometrica) either coincide with
the classical theorems just mentioned, or contain extra assumptions that are
alien to the Gauss-Markov or Aitken Theorem.
\end{abstract}

\section{Introduction}

\cite{Hansen_1.5, Hansen_2, Hansen_3} contain several assertions from which
he claims it would follow that the linearity condition can be dropped from
the Gauss-Markov Theorem or from the Aitken Theorem. We show that this
conclusion is unwarranted, as his assertions on which this conclusion rests
turn out to be only (intransparent) reformulations of the classical
Gauss-Markov or the classical Aitken Theorem, into which he has reintroduced
linearity through the backdoor, or contain extra assumptions alien to the
Gauss-Markov or Aitken Theorem.

The present paper is mainly pedagogical in nature. In particular, the
results will not come as a surprise to anyone well-versed in the theory of
linear models and familiar with basic concepts of statistical decision
theory, but -- given the confusion introduced by \cite{Hansen_1.5, Hansen_2,
Hansen_3} -- the paper will benefit the econometrics community.

One important upshot of the present paper is that one should \emph{not }%
follow Hansen's plea to drop the linearity condition in teaching the
Gauss-Markov Theorem or the Aitken Theorem. Depending on which formulation
of the Gauss-Markov Theorem one starts with (Theorem \ref{GM} or \ref%
{GMHansen} given below), dropping linearity from the formulation of that
theorem at best leads to a result equivalent to the usual Gauss-Markov
Theorem, and at worst leads to an incorrect result. The same goes for the
Aitken Theorem. Unfortunately, in heeding his own advice Hansen has included
an incorrect formulation of the Gauss-Markov Theorem in the August 2021
version of his forthcoming text-book (Theorem 4.4. in \cite{HABook})
available on his webpage for an extended time period.

\cite{Hansen_1.5} is the version accepted by Econometrica and which has been
available on Econometrica's webpage of forthcoming papers. \cite{Hansen_2}
is an updated version that corrects an incorrect result in \cite{Hansen_1.5}
(but otherwise is identical to the latter paper), and is available from
Hansen's webpage. \cite{Hansen_3} refers to the version finally published in
Econometrica, which contains several nontrivial changes relative to \cite%
{Hansen_1.5, Hansen_2} introduced into the paper at the proof-reading stage.
Because \cite{Hansen_1.5, Hansen_2} have been widely circulated and
discussed, and because \cite{Hansen_3} has been published in Econometrica,
there is a need to discuss all three versions. We shall start by first
discussing \cite{Hansen_1.5, Hansen_2} in Sections \ref{sect:GM} and \ref%
{Ait}, a discussion that has considerable bearings also on \cite{Hansen_3}.
We then move on to discuss the changes introduced into \cite{Hansen_3} at
the proof-reading stage and their ramifications in Section \ref{H3}. Section %
\ref{sect:iid} discusses the situation when one restricts attention to
independent identically distributed errors.

After the first version of this paper had been circulated, Stephen Portnoy
sent us a paper of his (\cite{Portnoy}) that has a result somewhat similar
to our Theorem \ref{class_trivial} with a different proof. For a discussion
see Section \ref{sect:GM}.

\section{The Framework\label{frame}}

As in \cite{Hansen_1.5, Hansen_2} we consider throughout the paper the
linear regression model 
\begin{equation}
Y=X\beta +e  \label{model}
\end{equation}%
where $Y$ is of dimension $n\times 1$ and $X$ is a (non-random) $n\times k$
design matrix with full column rank $k$ satisfying $1\leq k<n$.\footnote{%
We make the assumption $k<n$ in order to use exactly the same framework as
in Hansen's papers.} It is assumed that%
\begin{equation}
Ee=0  \label{expect}
\end{equation}%
and 
\begin{equation}
Eee^{\prime }=\sigma ^{2}\Sigma ,  \label{var}
\end{equation}%
where $\sigma ^{2}$, $0<\sigma ^{2}<\infty $, is unknown and $\Sigma $ is a
known symmetric and positive definite $n\times n$ matrix ($Ee^{\prime
}e<\infty $). This model implies a distribution $F$ for $Y$, which, for the
given $X$, depends on $\beta $ and the distribution of $e$, in particular on 
$\sigma ^{2}$ and $\Sigma $. Now define $\mathbf{F}_{2}(\Sigma )$ as the
class of all such distributions $F$ when $\beta $ varies through $\mathbb{R}%
^{k}$ and the distribution of $e$ varies through all distributions
compatible with (\ref{expect}) and (\ref{var}) for the given $\Sigma $ (and
arbitrary $\sigma ^{2}$, $0<\sigma ^{2}<\infty $). We furthermore introduce
the set $\mathbf{F}_{2}$ as the larger class where we also vary $\Sigma $
through the set of all symmetric and positive definite $n\times n$ matrices.
In other words,%
\begin{equation*}
\mathbf{F}_{2}=\tbigcup_{\Sigma }\mathbf{F}_{2}(\Sigma ),
\end{equation*}%
where the union is taken over all symmetric and positive definite $n\times n$
matrices.\footnote{%
Note that $\mathbf{F}_{2}(\Sigma _{1})\cap \mathbf{F}_{2}(\Sigma
_{2})=\emptyset $ iff $\Sigma _{1}$ and $\Sigma _{2}$ are not proportional.
And $\mathbf{F}_{2}(\Sigma _{1})=\mathbf{F}_{2}(\Sigma _{2})$ iff $\Sigma
_{1}$ and $\Sigma _{2}$ are proportional.} [Of course, $\mathbf{F}%
_{2}(\Sigma )$ as well as $\mathbf{F}_{2}$ also depend on the given $X$, but
this dependence is not shown in the notation.] The set $\mathbf{F}_{2}^{0}$
defined in \cite{Hansen_1.5, Hansen_2} is nothing else than $\mathbf{F}%
_{2}(I_{n})$, where $I_{n}$ denotes the $n\times n$ identity matrix.%
\footnote{%
Note that in \cite{Hansen_3} the symbol $\mathbf{F}_{2}^{0}$ is used to
denote a \emph{different} set of distributions; see Section \ref{H3} below.}
In the following $E_{F}$ ($Var_{F}$, respectively) will denote the
expectation (variance-covariance matrix, respectively) taken under the
distribution $F$. A word on notation: Given $F\in \mathbf{F}_{2}$, there is
a unique $\beta $, denoted by $\beta (F)$, and a unique $\sigma ^{2}\Sigma $%
, denoted by $(\sigma ^{2}\Sigma )(F)$, compatible with the distribution $F$.

\begin{remark}
\emph{(Ambiguity in the definition in \cite{Hansen_1.5, Hansen_2})\ }\cite%
{Hansen_1.5, Hansen_2} also define a set $\mathbf{F}_{2}$, unfortunately
somewhat ambiguously: Taking the first sentence mentioning his set $\mathbf{F%
}_{2}$ literally, his set would coincide with our $\mathbf{F}_{2}(\Sigma )$.
The two sentences following that sentence, however, intimate that his set $%
\mathbf{F}_{2}$ was intended to coincide with our set $\mathbf{F}_{2}$. This
is confirmed by an inspection of his proofs; furthermore, if one would
interpret his set $\mathbf{F}_{2}$ to mean our $\mathbf{F}_{2}(\Sigma )$,
then the relation $\mathbf{F}_{2}^{0}\subset \mathbf{F}_{2}$ given below (4)
in \cite{Hansen_1.5, Hansen_2} (which in our notation would become $\mathbf{F%
}_{2}(I_{n})\subset \mathbf{F}_{2}(\Sigma )$) could not hold (except for $%
\Sigma $ proportional to $I_{n}$). In the following we hence interpret
Hansen's set $\mathbf{F}_{2}$ to coincide with our definition of $\mathbf{F}%
_{2}$. In a remark further below we discuss what happens if one would adopt
the interpretation of Hansen's $\mathbf{F}_{2}$ as coinciding with our $%
\mathbf{F}_{2}(\Sigma )$.
\end{remark}

\section{The Gauss-Markov Case\label{sect:GM}}

To focus the discussion, we first treat the situation of a regression model
with homoskedastic and uncorrelated errors, i.e., we assume that in (\ref%
{var}) we have%
\begin{equation}
\Sigma =I_{n}.  \label{homo}
\end{equation}

Let $\hat{\beta}_{OLS}=(X^{\prime }X)^{-1}X^{\prime }Y$ denote the ordinary
least-squares estimator. The classical Gauss-Markov Theorem then reads as
follows. Recall that a linear estimator is of the form $AY$, were $A$ is a
(nonrandom) $k\times n$ matrix. Also recall that $\mathbf{F}_{2}^{0}=\mathbf{%
F}_{2}(I_{n})$.

\begin{theorem}
\label{GM}If $\hat{\beta}$ is a linear estimator that is unbiased under all $%
F\in \mathbf{F}_{2}^{0}$ (meaning that $E_{F}\hat{\beta}=\beta (F)$ for
every $F\in \mathbf{F}_{2}^{0}$), then 
\begin{equation*}
Var_{F}(\hat{\beta})\succeq Var_{F}(\hat{\beta}_{OLS})
\end{equation*}%
for every $F\in \mathbf{F}_{2}^{0}$. [Here $\succeq $ denotes Loewner order.]
\end{theorem}

The theorem can \emph{equivalently} be stated in the following more unusual
form, which is the form chosen by Hansen (see Theorem 1 in \cite{Hansen_1.5,
Hansen_2}).\footnote{\label{FN1}As formulated in \cite{Hansen_1.5, Hansen_2}%
, his\ Theorem 1 has $\sigma ^{2}(X^{\prime }X)^{-1}$ instead of $Var_{F}(%
\hat{\beta}_{OLS})$ on the r.h.s. of the inequality. Taken literally this
leaves $\sigma ^{2}$ unspecified. To obtain a mathematically well-defined
statement $\sigma ^{2}$ needs to be interpreted as $\sigma ^{2}(F)$, the
variance of the data under $F$, the distribution under which the
variance-covariance matrices of the estimators are computed.}

\begin{theorem}
\label{GMHansen}If $\hat{\beta}$ is a linear estimator that is unbiased
under all $F\in \mathbf{F}_{2}$ (meaning that $E_{F}\hat{\beta}=\beta (F)$
for every $F\in \mathbf{F}_{2}$), then 
\begin{equation}
Var_{F}(\hat{\beta})\succeq Var_{F}(\hat{\beta}_{OLS})  \label{varineq}
\end{equation}%
for every $F\in \mathbf{F}_{2}^{0}$.
\end{theorem}

In the latter theorem the unbiasedness is requested to hold over the \emph{%
larger} class $\mathbf{F}_{2}$ of distributions rather than only over $%
\mathbf{F}_{2}^{0}$. Of course, this is immaterial here and the two theorems
are equivalent, because the estimators are required to be \emph{linear} in
both theorems and thus their expectations depend only on the first moment of 
$Y$ and not on the second moments at all.\footnote{\label{FnA}For \emph{%
linear} estimators $\hat{\beta}$ the condition $E_{F}\hat{\beta}=\beta (F)$
for every $F\in \mathbf{F}_{2}$ is, in fact also equivalent to $E_{F}\hat{%
\beta}=\beta (F)$ for every $F\in \mathbf{G}$, whenever $\mathbf{G\subseteq }
$ $\mathbf{F}_{2}$ holds and $\{\beta (F):F\in \mathbf{G}\}$ contains a
basis of $\mathbb{R}^{k}$. This is obvious since any of these unbiasedness
conditions are equivalent to $AX=I_{k}$, where $A$ is the matrix
representing the linear estimator $\hat{\beta}$, i.e., $\hat{\beta}=AY$.
[For $\mathbf{G=F}_{2}(I_{n})$ we obtain the equivalence noted above in the
main text; a similar equivalence is obtained for $\mathbf{G=F}_{2}(\Sigma )$%
.] Furthermore, if $\mathbf{G}$ is chosen to correspond to all distributions
in $\mathbf{F}_{2}(I_{n})$ such that $e/\sigma $ follows a \emph{given}
distribution (\textquotedblleft parametric linear regression
model\textquotedblright ), the before noted equivalence applies, and we thus
can obtain a version of the Gauss-Markov Theorem for the parametric linear
regression model. A similar remark applies to the Aitken Theorem.} While the
difference in the unbiasedness conditions is immaterial in the preceding
theorems, it is worth pointing out that the unbiasedness condition as given
in Theorem \ref{GMHansen} requires that an estimator is not only unbiased in
the underlying model with uncorrelated and homoskedastic errors one is
studying, but also requires unbiasedness under correlated and/or
heteroskedastic errors (i.e., under structures that are `outside' of the
model that is being considered). Why one would want to impose such a
requirement when the underlying model has uncorrelated and homoskedastic
errors is at least debatable. However, we stress once more that in the
context of the preceding two theorems this does not matter due to the
assumed linearity of the estimators.

We next discuss what happens if one eliminates the linearity condition in
the two equivalent theorems. Dropping the linearity conditions leads to the
following assertions, which will turn out to be \emph{no longer} equivalent
to each other:

\bigskip

\textbf{Assertion 1:} If $\hat{\beta}$ is an estimator (i.e., a
Borel-measurable function of $Y$) that is unbiased under all $F\in \mathbf{F}%
_{2}^{0}$ (meaning that $E_{F}\hat{\beta}=\beta (F)$ for every $F\in \mathbf{%
F}_{2}^{0}$), then 
\begin{equation}
Var_{F}(\hat{\beta})\succeq Var_{F}(\hat{\beta}_{OLS})  \label{ineq_loew}
\end{equation}%
for every $F\in \mathbf{F}_{2}^{0}$.

\bigskip

\textbf{Assertion 2:} If $\hat{\beta}$ is an estimator that is unbiased
under all $F\in \mathbf{F}_{2}$ (meaning that $E_{F}\hat{\beta}=\beta (F)$
for every $F\in \mathbf{F}_{2}$), then 
\begin{equation*}
Var_{F}(\hat{\beta})\succeq Var_{F}(\hat{\beta}_{OLS})
\end{equation*}%
for every $F\in \mathbf{F}_{2}^{0}$.

\bigskip

Before discussing Assertions 1 and 2, we need to make a remark on the
interpretation of inequalities like (\ref{ineq_loew}).

\begin{remark}
\label{rem_Loew}(i) In Theorems \ref{GM} and \ref{GMHansen} the objects $%
Var_{F}(\hat{\beta})$ as well as $Var_{F}(\hat{\beta}_{OLS})$ are
well-defined as real matrices because all estimators considered are linear,
and hence $E_{F}(\parallel \hat{\beta}\parallel ^{2})<\infty $, $%
E_{F}(\parallel \hat{\beta}_{OLS}\parallel ^{2})<\infty $ holds for every $%
F\in \mathbf{F}_{2}^{0}$ where $\parallel .\parallel $ denotes the Euclidean
norm. In contrast, in Assertions 1 and 2 estimators $\hat{\beta}$ with $%
E_{F}(\parallel \hat{\beta}\parallel ^{2})=\infty $ for some $F\in \mathbf{F}%
_{2}^{0}$ are permissible. [Note that $E_{F}(\parallel \hat{\beta}\parallel
^{2})=\infty $ for some $F\in \mathbf{F}_{2}^{0}$ and $E_{F}(\parallel \hat{%
\beta}\parallel ^{2})<\infty $ for some other $F\in \mathbf{F}_{2}^{0}$ may
occur.] This necessitates some discussion how Assertions 1 and 2 are then to
be read. In the first version of this paper we unfortunately had glossed
over this issue, but an explicit discussion is warranted. For the subsequent
discussion note that in both assertions $E_{F}(\hat{\beta})$ is well-defined
and finite for every $F\in \mathbf{F}_{2}^{0}$ as a consequence of the
respective unbiasedness assumption (and because $\mathbf{F}_{2}^{0}\subseteq 
\mathbf{F}_{2}$).

(ii) In the \emph{scalar} case (i.e., $k=1$), there is no problem as the
object $Var_{F}(\hat{\beta})$ is well-defined for every $F\in \mathbf{F}%
_{2}^{0}$ as an element of the \emph{extended} real line, regardless of
whether $E_{F}(\parallel \hat{\beta}\parallel ^{2})<\infty $ or not. Hence,
inequality (\ref{ineq_loew}) always makes sense in case $k=1$.

(iii) For general $k$, in case the estimator $\hat{\beta}$ satisfies $%
E_{F}(\parallel \hat{\beta}\parallel ^{2})<\infty $ for a given $F\in 
\mathbf{F}_{2}^{0}$, the object $Var_{F}(\hat{\beta})$ is well-defined as a
real matrix. Note that the inequality (\ref{ineq_loew}) can then
equivalently be expressed as $Var_{F}(c^{\prime }\hat{\beta})\geq
Var_{F}(c^{\prime }\hat{\beta}_{OLS})$ for every $c\in \mathbb{R}^{k}$.

(iv) In the case $k>1$, the object $Var_{F}(\hat{\beta})$ is not
well-defined if $E_{F}(\parallel \hat{\beta}\parallel ^{2})=\infty $ ($F\in 
\mathbf{F}_{2}^{0}$), and hence it is not immediately clear how (\ref%
{ineq_loew}) should then be understood. However, the inequalities $%
Var_{F}(c^{\prime }\hat{\beta})\geq Var_{F}(c^{\prime }\hat{\beta}_{OLS})$
for every $c\in \mathbb{R}^{k}$ still make sense in view of (ii) above. We
hence may and will interpret (\ref{ineq_loew}) (with $F\in \mathbf{F}%
_{2}^{0} $) as a \emph{symbolic shorthand notation} for $Var_{F}(c^{\prime }%
\hat{\beta})\geq Var_{F}(c^{\prime }\hat{\beta}_{OLS})$ for every $c\in 
\mathbb{R}^{k}$ (which works both in the case $E_{F}(\parallel \hat{\beta}%
\parallel ^{2})<\infty $ and in the case $E_{F}(\parallel \hat{\beta}%
\parallel ^{2})=\infty $). We have chosen to write inequality (\ref%
{ineq_loew}) as given (abusing notation), rather than the more conventional
and more precise $Var_{F}(c^{\prime }\hat{\beta})\geq Var_{F}(c^{\prime }%
\hat{\beta}_{OLS})$ for every $c\in \mathbb{R}^{k}$, in order for our
discussion to be easily comparable with the presentation in Hansen's papers;
his papers are silent on this issue. The same convention applies mutatis
mutandis to similar statements such as, e.g., Assertions 3 and 4, etc.

(v) The above discussion would become moot, if one would introduce the \emph{%
extra} assumption $E_{F}(\parallel \hat{\beta}\parallel ^{2})<\infty $ for
every $F\in \mathbf{F}_{2}^{0}$ into Assertions 1 and 2. However, such an
additional assumption, which has little justification, would (potentially)
narrow down the class of estimators competing with $\hat{\beta}_{OLS}$. As
we shall see later on, such an extra assumption actually would have no
effect on Assertion 2 (and thus on the corresponding theorems in Hansen's
papers) at all in view of our Theorem \ref{class_trivial}. The effect it
would have on Assertion 1 (and some other results) is discussed in Appendix %
\ref{App_B}.
\end{remark}

We now turn to discussing Assertions 1 and 2. Not unexpectedly, Assertion 1
is \emph{incorrect} in general.\footnote{%
I.e., there exist design matrices $X$ such that the assertion is false.}
This is known. For the benefit of the reader we provide some counterexamples
and attending discussion in Appendix \ref{App_A}. In particular, we see that
in the classical Gauss-Markov Theorem \emph{as it is usually formulated}
(i.e., in Theorem \ref{GM}) one can \emph{not }eliminate the linearity
condition in general!

Concerning Assertion 2, note that it coincides with Theorem 5 in \cite%
{Hansen_1.5, Hansen_2} (his `modern Gauss-Markov Theorem').\footnote{%
The same caveat as expressed in Footnote \ref{FN1} also applies to the
formulation of Theorem 5 in \cite{Hansen_1.5, Hansen_2}.} Obvious questions
now are (i) whether Assertion 2 (i.e., Theorem 5 in \cite{Hansen_1.5,
Hansen_2}) is correct, and (ii) if so, what is the reason for Assertion 2 to
be correct while Assertion 1 is incorrect in general although in both
assertions the linearity condition has been dropped. The answer to the
latter question lies in the fact that Assertion 2 is requiring a \emph{%
stricter} unbiasedness condition, namely unbiasedness over $\mathbf{F}_{2}$
rather than only unbiasedness over $\mathbf{F}_{2}^{0}$. While the two
unbiasedness conditions effectively coincide for \emph{linear} estimators as
discussed before, this is no longer the case once we leave the realm of
linear estimators. Hence, the (potential) correctness of Assertion 2 (i.e.,
of Theorem 5 in \cite{Hansen_1.5, Hansen_2}) must crucially rest on imposing
the stricter unbiasedness condition, which not only requires unbiasedness
under the model considered (regression with homoskedastic and uncorrelated
errors), but oddly also under structures `outside' of the maintained model
(namely under heteroskedastic and/or correlated errors). Note that the class
of competitors to $\hat{\beta}_{OLS}$ figuring in Assertion 1 is, in
general, \emph{larger} than the class of competitors appearing in Assertion
2. Nevertheless, \cite{Hansen_1.5, Hansen_2} (and also \cite{Hansen_3}) are
quiet on the use of this stricter unbiasedness condition.

Having understood what distinguishes Assertion 2 (i.e., Theorem 5 in \cite%
{Hansen_1.5, Hansen_2}) from Assertion 1, the question remains whether the
former is indeed correct, and if so, what its scope is, i.e., how much
larger than the class of linear (unbiased) estimators the class of
estimators covered by Assertion 2 (i.e., by Theorem 5 in \cite{Hansen_1.5,
Hansen_2}) is. We answer this now: As we shall show in the subsequent
theorem, the only estimators $\hat{\beta}$ satisfying the unbiasedness
condition of Assertion 2 (i.e., of Theorem 5 in \cite{Hansen_1.5, Hansen_2})
are \emph{linear} estimators! \emph{In other words, Theorem 5 in \cite%
{Hansen_1.5, Hansen_2} (i.e., his `modern Gauss-Markov Theorem') is nothing
else than the good old(fashioned) Gauss-Markov Theorem (i.e., Theorem \ref%
{GM} above), just stated in a somewhat unusual and intransparent way!}%
\footnote{%
Recall from before that for linear estimators the unbiasedness conditions in
Theorems \ref{GM} and \ref{GMHansen} are equivalent.} [While the word
`linear' does not appear in the \emph{formulation} of Theorem 5 in \cite%
{Hansen_1.5, Hansen_2}, linearity of the estimators is introduced indirectly
through a backdoor provided by the stricter unbiasedness condition.] While
Theorem 5 in \cite{Hansen_1.5, Hansen_2} thus turns out to be correct, it is
certainly not new!\footnote{%
We have not checked whether the \emph{proofs} in \cite{Hansen_1.5, Hansen_2}
are correct or not.} Theorem 6 in \cite{Hansen_2} is a special case of his
Theorem 5 for the location model, and thus is also not new; in contrast,
Theorem 6 in \cite{Hansen_1.5} is a special case of Assertion 1. Example \ref%
{EX1} in Appendix \ref{App_A} shows that this theorem is false. What has
been said so far also serves as a reminder that one has to be careful with
statements such as \textquotedblleft best unbiased equals best linear
unbiased\textquotedblright . While this statement is incorrect in the
context of Assertion 1 in general, it is trivially correct in the context of
Assertion 2 (i.e., of Theorem 5 in \cite{Hansen_1.5, Hansen_2}) as a
consequence of the subsequent Theorem \ref{class_trivial}.

An upshot of the preceding discussion is that -- despite a plea to the
contrary in \cite{Hansen_1.5, Hansen_2, Hansen_3} -- one should \textbf{not}
drop `linearity' from the pedagogy of the Gauss-Markov Theorem. There is
nothing to gain and a lot to lose: It will lead to an incorrect assertion,
if one starts from the usual formulation of the classical Gauss-Markov
Theorem (i.e., from Theorem \ref{GM});\ otherwise (i.e., if one starts from
Theorem \ref{GMHansen}), it will lead to a correct, but rather
intransparent, assertion that is in fact equivalent to the classical
Gauss-Markov Theorem. Unfortunately, Hansen has fallen victim to his own
advice as the Gauss-Markov Theorem (Theorem 4.4) given in the August 2021
version of his forthcoming text-book \cite{HABook} is incorrect in general
(as it coincides with Assertion 1).

We now provide the theorem alluded to above. After the first version of this
paper had been circulated, we learned about \cite{Portnoy}, which
establishes a related result using different arguments than the ones we use;
for more discussion see Remark \ref{rem_P} further below.

\begin{theorem}
\label{class_trivial}If $\hat{\beta}$ is an estimator (i.e., a
Borel-measurable function of $Y$) that is unbiased under all $F\in \mathbf{F}%
_{2}$ (meaning that $E_{F}\hat{\beta}=\beta (F)$ for every $F\in \mathbf{F}%
_{2}$), then $\hat{\beta}$ is a linear estimator (i.e., $\hat{\beta}=AY$ for
some $k\times n$ matrix $A$).\footnote{%
By unbiasedness, such an $A$ must then also satisfy $AX=I_{k}$.}
\end{theorem}

We give a first "proof" based on Theorem 4.3 in \cite{koopmann} (also
reported as Theorem 2.1 in \cite{Gnotetal}), but see the discussion
immediately following this "proof" for a caveat.\footnote{%
Curiously, the result by \cite{koopmann} in question is actually mentioned
in Section 1 of \cite{Hansen_1.5, Hansen_2, Hansen_3}.}

\textbf{A first "proof":} The unbiasedness assumption of the theorem
obviously translates into%
\begin{equation}
E_{F}\hat{\beta}=\beta (F)\text{ \ for every \ }F\in \mathbf{F}_{2}(\Sigma ),
\label{koop}
\end{equation}%
for \emph{every} symmetric and positive definite $\Sigma $ of dimension $%
n\times n$; specializing to the case $\Sigma =I_{n}$, we, in particular,
obtain\footnote{%
Instead of $I_{n}$ we could have chosen any other symmetric and positive
definite $n\times n$ matrix $\Sigma _{0}$ instead.}%
\begin{equation}
E_{F}\hat{\beta}=\beta (F)\text{ \ for every \ }F\in \mathbf{F}_{2}(I_{n}).
\label{koop2}
\end{equation}%
Condition (\ref{koop2}), together with Theorem 4.3 in \cite{koopmann} (see
also Theorem 2.1 in \cite{Gnotetal}\footnote{%
Note that $X^{-}$ in that reference runs through all possible $g$-inverses
of $X$.}$^{\text{,}}$\footnote{\label{fn_koop}\cite{Gnotetal} assume $\sigma
^{2}>0$ whereas \cite{koopmann} allows also $\sigma ^{2}=0$. However, both
theorems are equivalent as unbiasedness under every $F\in \mathbf{F}%
_{2}(I_{n})$ also implies unbiasedness under the point distributions at $%
X\beta $ (i.e., the distributions corresponding to $\sigma ^{2}=0$). This is
easily seen by considering those distributions in $\mathbf{F}_{2}(I_{n})$
that correspond to $X\beta +e$ with the components of $e$ being independent
identically distributed according to $\varepsilon _{m}(\delta _{-1}+\delta
_{1})/2+(1-\varepsilon _{m})\delta _{0}$. Here $\varepsilon _{m}$, $%
0<\varepsilon _{m}<1$, converges to zero for $m\rightarrow \infty $ and $%
\delta _{x}$ denotes point mass at $x\in \mathbb{R}$. A similar argument
applies in the case of $\mathbf{F}_{2}(\Sigma )$.}), implies that $\hat{\beta%
}$ is of the form%
\begin{equation}
\hat{\beta}=A^{0}Y+(Y^{\prime }H_{1}^{0}Y,\ldots ,Y^{\prime
}H_{k}^{0}Y)^{\prime },  \label{quadr}
\end{equation}%
where $A^{0}$ satisfies $A^{0}X=I_{k}$ and $H_{i}^{0}$ are matrices
satisfying $\limfunc{tr}(H_{i}^{0})=0$ and $X^{\prime }H_{i}^{0}X=0$ for $%
i=1,\ldots ,k$. It is easy to see that we may without loss of generality
assume that the matrices $H_{i}^{0}$ are symmetric (otherwise replace $%
H_{i}^{0}$ by $(H_{i}^{0}+H_{i}^{0\prime })/2$). Inserting (\ref{quadr})
into (\ref{koop}) yields%
\begin{equation*}
E_{F}\left( A^{0}Y+(Y^{\prime }H_{1}^{0}Y,\ldots ,Y^{\prime
}H_{k}^{0}Y)^{\prime }\right) =\beta (F)\text{ \ for every \ }F\in \mathbf{F}%
_{2}(\Sigma ),
\end{equation*}%
and this has to hold for \emph{every} symmetric and positive definite $%
\Sigma $. Standard calculations involving the trace operator and division by 
$\sigma ^{2}$ now give%
\begin{equation}
(\limfunc{tr}(H_{1}^{0}\Sigma ),\ldots ,\limfunc{tr}(H_{k}^{0}\Sigma
))^{\prime }=0\text{ \ for every symmetric and positive definite }\Sigma .
\label{intermed}
\end{equation}%
For every $j=1,\ldots ,n$, choose now a sequence of symmetric and positive
definite matrices $\Sigma _{m}^{(j)}$ (each of dimension $n\times n)$ that
converges to $e_{j}(n)e_{j}(n)^{\prime }$ as $m\rightarrow \infty $, where $%
e_{j}(n)$ denotes the $j$-th standard basis vector in $\mathbb{R}^{n}$ (such
sequences obviously exist). Plugging this sequence into (\ref{intermed}),
letting $m$ go to infinity, and exploiting properties of the trace-operator,
we obtain%
\begin{equation*}
(e_{j}(n)^{\prime }H_{1}^{0}e_{j}(n),\ldots ,e_{j}(n)^{\prime
}H_{k}^{0}e_{j}(n))^{\prime }=0\text{ \ for every }j=1,\ldots ,n\text{.}
\end{equation*}%
In other words, all the diagonal elements of $H_{i}^{0}$ are zero for every $%
i=1,\ldots ,k$. Next, for every $j,l=1,\ldots ,n$, $j\neq l$, choose a
sequence of symmetric and positive definite matrices $\Sigma _{m}^{\{j,l\}}$
(each of dimension $n\times n)$ that converges to $%
(e_{j}(n)+e_{l}(n))(e_{j}(n)+e_{l}(n))^{\prime }$ as $m\rightarrow \infty $
(such sequences obviously exist). Then exactly the same argument as before
delivers%
\begin{equation*}
((e_{j}(n)+e_{l}(n))^{\prime }H_{1}^{0}(e_{j}(n)+e_{l}(n)),\ldots
,(e_{j}(n)+e_{l}(n))^{\prime }H_{k}^{0}(e_{j}(n)+e_{l}(n)))^{\prime }=0\text{
\ for every }j\neq l\text{.}
\end{equation*}%
Recall that the matrices $H_{i}^{0}$ are symmetric. Together with the
already established fact that the diagonal elements are all zero, we obtain
that also all the off-diagonal elements in any of the matrices $H_{i}^{0}$
are zero; i.e., $H_{i}^{0}=0$ for every $i=1,\ldots ,k$. This completes the
proof.\footnote{%
A slightly different version of the first "proof" can be obtained as
follows. Theorem 4.3 in \cite{koopmann} (together with Footnote \ref{fn_koop}%
) shows for every given (fixed) $\Sigma $ that any $\hat{\beta}$ satisfying (%
\ref{koop}) is of the form $AY+(Y^{\prime }H_{1}Y,\ldots ,Y^{\prime
}H_{k}Y)^{\prime }$ where $AX=I_{k}$, the $H_{i}$'s satisfy $\limfunc{tr}%
(H_{i}\Sigma )=0$, and $X^{\prime }H_{i}X=0$ for $i=1,\ldots ,k$. Again it
is easy to see that we may assume the matrices $H_{i}$ to be symmetric. Note
that the matrices $A$ and $H_{i}$ flowing from Theorem 4.3 in \cite{koopmann}
in principle could depend on $\Sigma $. The following argument shows that
this is, however, not the case (after symmetrization of the $H_{i}$'s) in
the present situation: If $\hat{\beta}$ had two distinct linear-quadratic
representations with symmetric $H_{i}$'s, then the difference of these two
representations would be a vector of multivariate polynomials (at least one
of which is nontrivial) that would have to vanish everywhere, which is
impossible since the zero-set of a nontrivial multivariate polynomial is a
Lebesgue null-set. Given now the independence (from $\Sigma $) of the
matrices $H_{i}$, one can then exploit the before mentioned relations $%
\limfunc{tr}(H_{i}\Sigma )=0$ in the same way as is done following (\ref%
{intermed}) in the main text.} $\blacksquare $

Theorem 4.3 in \cite{koopmann} is proved by reducing it to Theorem 3.1 (via
Theorems 3.2, 4.1, and 4.2) in the same reference. Unfortunately, a full
proof of Theorem 3.1 is not provided in \cite{koopmann}, only a very rough
outline is given. Thus the status of Theorem 4.3 in \cite{koopmann} is not
entirely clear. For this reason we next give a direct proof of our Theorem %
\ref{class_trivial} which does not rely on any result in \cite{koopmann}.%
\footnote{%
Alternatively, one could try to provide a complete proof of the result in 
\cite{koopmann}. We have not pursued this, but have chosen the route via a
direct proof of our Theorem \ref{class_trivial}.}

\textbf{A direct proof:} It suffices to establish $\hat{\beta}(y+z)=\hat{%
\beta}(y)+\hat{\beta}(z)$ as well as $\hat{\beta}(cz)=c\hat{\beta}(z)$ for
every $y$ and $z$ in $\mathbb{R}^{n}$ and every $c\in \mathbb{R}$. For every 
$m\in \mathbb{N}$ with $m\geq 2$, every $V=(v_{1},\ldots ,v_{m})\in \mathbb{R%
}^{n\times m}$ and $\alpha \in (0,1)^{m}$ such that $\sum_{i=1}^{m}\alpha
_{i}=1$, define a probability measure (distribution) via%
\begin{equation*}
\mu _{V,\alpha }:=\sum_{i=1}^{m}\alpha _{i}\delta _{v_{i}},
\end{equation*}%
where $\delta _{z}$ denotes unit point mass at $z\in \mathbb{R}^{n}$. The
expectation of $\mu _{V,\alpha }$ equals $V\alpha $, and its
variance-covariance matrix equals $V\limfunc{diag}(\alpha )V^{\prime
}-(V\alpha )(V\alpha )^{\prime }$. Denote the expectation operator w.r.t.~$%
\mu _{V,\alpha }$ by $E_{V,\alpha }$. Note that in case $V\alpha =0$ and $%
\limfunc{rank}(V)=n$ the measure $\mu _{V,\alpha }$ has expectation zero and
a positive definite variance-covariance matrix; thus, $\mu _{V,\alpha }$
corresponds to an $F\in \mathbf{F}_{2}$ which has $\beta (F)=0$. From the
unbiasedness assumption imposed on $\hat{\beta}$ we obtain that%
\begin{equation}
V\alpha =0\text{ and\ }\limfunc{rank}(V)=n\text{ \ implies \ }0=E_{V,\alpha
}(\hat{\beta})=\sum_{i=1}^{m}\alpha _{i}\hat{\beta}(v_{i}).  \label{eqn:o1}
\end{equation}

\textbf{Step 1:} Fix $z\in \mathbb{R}^{n}$ and define $\alpha
^{(1)}=2^{-1}(n^{-1},\ldots ,n^{-1})^{\prime }\in \mathbb{R}^{2n}$, $\alpha
^{(2)}=2^{-1}((n+1)^{-1},\ldots ,(n+1)^{-1})^{\prime }\in \mathbb{R}%
^{2(n+1)} $, $V_{1}=(I_{n},-I_{n})$ and $V_{2}=(I_{n},-I_{n},z,-z)$. Clearly 
$V_{1}\alpha ^{(1)}=V_{2}\alpha ^{(2)}=0$ and $\limfunc{rank}(V_{1})=%
\limfunc{rank}(V_{2})=n$. Furthermore, 
\begin{equation}
\mu _{V_{2},\alpha ^{(2)}}=\frac{n}{n+1}\mu _{V_{1},\alpha ^{(1)}}+\frac{1}{%
2(n+1)}(\delta _{z}+\delta _{-z}),  \label{eqn:o2}
\end{equation}%
which implies%
\begin{equation*}
E_{V_{2},\alpha ^{(2)}}(\hat{\beta})=\frac{n}{n+1}E_{V_{1},\alpha ^{(1)}}(%
\hat{\beta})+\frac{1}{2(n+1)}(\hat{\beta}(z)+\hat{\beta}(-z)).
\end{equation*}%
Applying (\ref{eqn:o1}) to $E_{V_{2},\alpha ^{(2)}}(\hat{\beta})$ and $%
E_{V_{1},\alpha ^{(1)}}(\hat{\beta})$ now yields $0=\hat{\beta}(z)+\hat{\beta%
}(-z)$, i.e., we have shown that 
\begin{equation}
\hat{\beta}(-z)=-\hat{\beta}(z)\text{ for every }z\in \mathbb{R}^{n},
\label{eqn:unev}
\end{equation}%
in particular $\hat{\beta}(0)=0$ follows.

\textbf{Step 2:} Let $y$ and $z$ be elements of $\mathbb{R}^{n}$. Define the
matrix 
\begin{equation*}
A(y,z)=((y_{1}+z_{1})e_{1}(n),\ldots ,(y_{n}+z_{n})e_{n}(n)),
\end{equation*}%
where $e_{i}(n)$ denotes the $i$-th standard basis vector in $\mathbb{R}^{n}$%
, and set 
\begin{equation*}
V=(A(y,z),-y,-z,I_{n},-I_{n})\quad \text{ and }\quad \alpha
=(3n+2)^{-1}(1,\ldots ,1)^{\prime }\in \mathbb{R}^{3n+2}.
\end{equation*}%
Then, we obtain $V\alpha =0$ and $\limfunc{rank}(V)=n$. Using (\ref{eqn:o1})
and (\ref{eqn:unev}) it follows that 
\begin{equation*}
0=\sum_{i=1}^{n}\hat{\beta}((y_{i}+z_{i})e_{i}(n))+\hat{\beta}(-y)+\hat{\beta%
}(-z),
\end{equation*}%
which by (\ref{eqn:unev}) is equivalent to 
\begin{equation}
\hat{\beta}(y)+\hat{\beta}(z)=\sum_{i=1}^{n}\hat{\beta}%
((y_{i}+z_{i})e_{i}(n)).  \label{eqn:rel}
\end{equation}%
Using (\ref{eqn:rel}) with $y$ replaced by $y+z$ and $z$ replaced by $0$
yields 
\begin{equation*}
\hat{\beta}(y+z)+\hat{\beta}(0)=\sum_{i=1}^{n}\hat{\beta}%
((y_{i}+z_{i})e_{i}(n)).
\end{equation*}%
Since $\hat{\beta}(0)=0$ as shown before, we obtain 
\begin{equation}
\hat{\beta}(y)+\hat{\beta}(z)=\hat{\beta}(y+z)\quad \text{ for every }y\text{
and }z\text{ in }\mathbb{R}^{n}.  \label{Cauchy}
\end{equation}%
That is, we have shown that $\hat{\beta}$ is additive, i.e., is a group
homomorphism between the additive groups $\mathbb{R}^{n}$ and $\mathbb{R}%
^{k} $. By assumption it is also Borel-measurable. It then follows by a
result due to Banach and Pettis (e.g., Theorem 2.2 in \cite{rosendal_2009})
that $\hat{\beta}$ is also continuous. Homogeneity of $\hat{\beta}$ now
follows from a standard argument, dating back to Cauchy, so that $\hat{\beta}
$ is in fact linear. We give the details for the convenience of the reader:
Relation (\ref{Cauchy}) (which contains (\ref{eqn:unev}) as a special case)
implies $\hat{\beta}(lz)=l\hat{\beta}(z)$ for every integer $l$. Replacing $%
z $ by $z/l $ ($l\neq 0$) in the latter relation gives $\hat{\beta}(z)/l=%
\hat{\beta}(z/l) $ for integer $l\neq 0$. It immediately follows that $\hat{%
\beta}(pz/q)=(p/q)\hat{\beta}(z)$ for every pair of integers $p$ and $q$ ($%
q\neq 0$). Let $c\in \mathbb{R}$ be arbitrary. Choose a sequence of rational
numbers $c_{s}$ that converges to $c$. Then by continuity of $\hat{\beta}$%
\begin{equation*}
\hat{\beta}(cz)=\lim_{s\rightarrow \infty }\hat{\beta}(c_{s}z)=\lim_{s%
\rightarrow \infty }\left( c_{s}\hat{\beta}(z)\right) =\left(
\lim_{s\rightarrow \infty }c_{s}\right) \hat{\beta}(z)=c\hat{\beta}(z).
\end{equation*}%
This concludes the proof. $\blacksquare $

\begin{remark}
\label{rem_prev}Inspection of the direct proof above shows that it does not
make use of the full force of the unbiasedness condition ($E_{F}\hat{\beta}%
=\beta (F)$ for every $F\in \mathbf{F}_{2}$), but only exploits unbiasedness
for certain strategically chosen discrete distributions $F$, each with
finite support and satisfying $\beta (F)=0$.
\end{remark}

\begin{remark}
\label{rem_P}(i) \cite{Portnoy} uses a somewhat weaker unbiasedness
condition than the one used in our Theorem \ref{class_trivial} (but see
Remark \ref{rem_prev}), and then establishes only Lebesgue almost everywhere
linearity of the estimators rather than linearity. This is an important
distinction for the following reason: The results in \cite{Hansen_1.5,
Hansen_2, Hansen_3} allow also for discrete distributions. For such
distributions positive probability mass can fall into the exceptional
Lebesgue null set, showing that any attempt to enforce linearity by
appropriately redefining the estimator on the exceptional null set will in
general not preserve the statistical properties of the estimator. In
particular, the claim in Comment (a) in Section 3 of \cite{Portnoy} that his
result \textquotedblleft implies Hansen's result\textquotedblright\ is not
warranted. Furthermore, at several instances in the discussion in \cite%
{Portnoy} linearity is incorrectly claimed although only linearity Lebesgue
almost everywhere is actually established in his paper.

(ii) \cite{Portnoy} emphasizes in his introduction as well as in Comment (b)
in his Section 3 that his result allows for distributions that have no
finite second moment. The following comment seems to be in order: The proof
in \cite{Portnoy} relies on requiring unbiasedness over a certain class $%
\mathbf{P}$, say, of distributions which have compact support, and thus have
finite moments of all orders. Trivially, then Portnoy's result holds a
fortiori if one requires unbiasedness to hold over a larger class $\mathbf{P}%
^{\ast }\mathbf{\supseteq P}$ of distributions, where $\mathbf{P}^{\ast }$
may contain also distributions that only have a finite first moment, but no
finite second moment. The direct proof of our Theorem 3.4 effectively relies
only on unbiasedness over a family of discrete distributions, each having
finite support (cf.~Remark \ref{rem_prev}). Again, then our linearity result
trivially holds a fortiori if unbiasedness is required over any class of
distributions containing the before mentioned family of discrete
distributions. Of course, such a class may then also contain distributions
that only have a finite first, but no finite second moment.
\end{remark}

\begin{remark}
\label{rem:amb}\emph{(Ambiguity in the definition in \cite{Hansen_1.5,
Hansen_2} continued)} If Hansen's $\mathbf{F}_{2}$ would be interpreted as
coinciding with our $\mathbf{F}_{2}(\Sigma )$ (here with $\Sigma =I_{n}$
because of (\ref{homo})) then the formulations of Theorems \ref{GM} and \ref%
{GMHansen} as well as the formulations of Assertions 1 and 2 would coincide.
In particular, with such an interpretation of Hansen's $\mathbf{F}_{2}$ his
Theorem 5 in \cite{Hansen_1.5, Hansen_2} would be false.
\end{remark}

\section{The Aitken Case\label{Ait}}

In this section we drop the assumption (\ref{homo}), i.e., $\Sigma $ in (\ref%
{var}) need not be the identity matrix. We make a preparatory remark:
Similarly to observations made in Section \ref{sect:GM} (see Footnote \ref%
{FN1}), the rendition of Aitken's Theorem (for linear estimators) as given
in Theorem 3 in \cite{Hansen_1.5, Hansen_2} needs some interpretation to
convert it into a mathematically well-defined statement: The product $\sigma
^{2}\Sigma $, on which the r.h.s.~of the inequality in that theorem depends
(note that $\sigma ^{2}$ and $\Sigma $ enter the expression only via the
product), is unspecified, and needs to be interpreted as $(\sigma ^{2}\Sigma
)(F)$, the variance-covariance matrix of the data under the relevant $F$
w.r.t.~which the variance-covariances in this inequality are taken. The same
comment applies to Theorem 4 in \cite{Hansen_1.5, Hansen_2}.

Aitken's Theorem as usually given in the literature reads as follows. Let $%
\hat{\beta}_{GLS}=\hat{\beta}_{GLS}(\Sigma )=(X^{\prime }\Sigma
^{-1}X)^{-1}X^{\prime }\Sigma ^{-1}Y$ denote the generalized least-squares
estimator using the known matrix $\Sigma $. Linear estimators are of the
form $\hat{\beta}=AY$ where $A$ is a (nonrandom) $k\times n$ matrix.

\begin{theorem}
\label{Aitken}Let $\Sigma $ be an arbitrary known symmetric and positive
definite $n\times n$ matrix. If $\hat{\beta}$ is a linear estimator that is
unbiased under all $F\in \mathbf{F}_{2}(\Sigma )$ (meaning that $E_{F}\hat{%
\beta}=\beta (F)$ for every $F\in \mathbf{F}_{2}(\Sigma )$), then 
\begin{equation*}
Var_{F}(\hat{\beta})\succeq Var_{F}(\hat{\beta}_{GLS})
\end{equation*}%
for every $F\in \mathbf{F}_{2}(\Sigma )$.
\end{theorem}

Similar as in Section \ref{sect:GM}, due to linearity of the estimators, an 
\emph{equivalent} version of the theorem is obtained if the unbiasedness
requirement is extended to all of $\mathbf{F}_{2}$.\footnote{%
Cf.~Footnote~\ref{FnA}.} This is precisely what happens in Theorem 3 in \cite%
{Hansen_1.5, Hansen_2}, his rendition of the Aitken Theorem (for linear
estimators). Note that the subsequent theorem is obviously equivalent to
Theorem 3 in \cite{Hansen_1.5, Hansen_2} and perhaps is more transparent.
[To see the equivalence, note that the all-quantor over $\Sigma $ in Theorem %
\ref{Aitken_variant} can be "absorbed" by replacing $\mathbf{F}_{2}(\Sigma )$
in that theorem with $\mathbf{F}_{2}$, provided the quantity $\sigma
^{2}\Sigma $ appearing in the expression $Var_{F}(\hat{\beta}_{GLS})=\sigma
^{2}(X^{\prime }\Sigma ^{-1}X)^{-1}=(X^{\prime }(\sigma ^{2}\Sigma
)^{-1}X)^{-1}$ in (\ref{vv}) below is understood as $(\sigma ^{2}\Sigma )(F)$%
, as is necessary anyway for Theorem 3 in \cite{Hansen_1.5, Hansen_2} to
formally make sense as noted earlier.]

\begin{theorem}
\label{Aitken_variant}Let $\Sigma $ be an arbitrary known symmetric and
positive definite $n\times n$ matrix. If $\hat{\beta}$ is a linear estimator
that is unbiased under all $F\in \mathbf{F}_{2}$ (meaning that $E_{F}\hat{%
\beta}=\beta (F)$ for every $F\in \mathbf{F}_{2}$), then 
\begin{equation}
Var_{F}(\hat{\beta})\succeq Var_{F}(\hat{\beta}_{GLS})  \label{vv}
\end{equation}%
for every $F\in \mathbf{F}_{2}(\Sigma )$.
\end{theorem}

Dropping linearity in both theorems now leads to two assertions.

\bigskip

\textbf{Assertion 3: }Let $\Sigma $ be an arbitrary known symmetric and
positive definite $n\times n$ matrix. If $\hat{\beta}$ is an estimator that
is unbiased under all $F\in \mathbf{F}_{2}(\Sigma )$ (meaning that $E_{F}%
\hat{\beta}=\beta (F)$ for every $F\in \mathbf{F}_{2}(\Sigma )$), then 
\begin{equation*}
Var_{F}(\hat{\beta})\succeq Var_{F}(\hat{\beta}_{GLS})
\end{equation*}%
for every $F\in \mathbf{F}_{2}(\Sigma )$.

\bigskip

\textbf{Assertion 4: }Let $\Sigma $ be an arbitrary known symmetric and
positive definite $n\times n$ matrix. If $\hat{\beta}$ is an estimator that
is unbiased under all $F\in \mathbf{F}_{2}$ (meaning that $E_{F}\hat{\beta}%
=\beta (F)$ for every $F\in \mathbf{F}_{2}$), then 
\begin{equation*}
Var_{F}(\hat{\beta})\succeq Var_{F}(\hat{\beta}_{GLS})
\end{equation*}%
for every $F\in \mathbf{F}_{2}(\Sigma )$.

\bigskip

Assertion 3 is again incorrect in general for reasons similar to the ones
given for Assertion 1 in the previous section, cf.~Appendix \ref{App_A}.
Assertion 4 is equivalent to Theorem 4 in \cite{Hansen_1.5, Hansen_2} (to
which we shall refer as his `modern Aitken Theorem'); this is seen in the
same way as the equivalence of Theorem \ref{Aitken_variant} above with
Theorem 3 in \cite{Hansen_1.5, Hansen_2}. Assertion 4 is indeed correct, but
again not new, as the class of estimators figuring in Assertion 4 consists
only of linear estimators as a consequence of Theorem \ref{class_trivial}
above.\footnote{%
Adding the extra condition $E_{F}(\parallel \hat{\beta}\parallel
^{2})<\infty $ for every $F\in \mathbf{F}_{2}(\Sigma )$ would have no effect
on Assertion 4 in view of our Theorem \ref{class_trivial}. The effect this
extra condition would have on Assertion 3 is discussed in Appendix \ref%
{App_B}.} Furthermore, a comment like Remark \ref{rem:amb} also applies
here. We conclude this section by noting that the rendition of Aitken's
Theorem in the text-book \cite{HABook} (Theorem 4.5) is ambiguously
formulated, making it difficult to decide whether it coincides with the
(incorrect) Assertion 3 or with Assertion 4, which is (trivially) correct.

\section{The Results in \protect\cite{Hansen_3} \label{H3}}

In \cite{Hansen_3} the same model given by (\ref{model}), (\ref{expect}),
and (\ref{var}) as in \cite{Hansen_1.5, Hansen_2} is considered and $\mathbf{%
F}_{2}$ is defined in the same manner.\footnote{%
The assumption in \cite{Hansen_1.5, Hansen_2} that $\Sigma $ is known and
positive definite and that $\sigma ^{2}$ is positive has been dropped in 
\cite{Hansen_3}. Nevertheless positive definiteness of $\Sigma $ as well as $%
\sigma ^{2}>0$ are frequently used in \cite{Hansen_3} (e.g., inverses of $%
\Sigma $ are taken; the proof of Theorem 4 makes use of both properties). We
hence will continue to assume $\sigma ^{2}>0$ and positive definiteness of $%
\Sigma $ in our discussion. We furthermore note that the ambiguity in the
definition of $\mathbf{F}_{2}$ in \cite{Hansen_1.5, Hansen_2} is now being
avoided in \cite{Hansen_3} as $\Sigma $ is no longer assumed to be known. Of
course, then $\sigma ^{2}$ and $\Sigma $ are no longer identifiable.} A set $%
\mathbf{F}_{2}^{\ast }$ representing the subset of $\mathbf{F}_{2}$
corresponding to \emph{independent} errors $e_{1},\ldots ,e_{n}$ is also
defined; here $e_{i}$ denotes the $i$-th component of the error vector $e$.
Furthermore, the subset of $\mathbf{F}_{2}^{\ast }$ corresponding to \emph{%
independent homoskedastic} errors is denoted by $\mathbf{F}_{2}^{0}$ in \cite%
{Hansen_3}. It should be noted that this set is \emph{not} the same as the
set $\mathbf{F}_{2}^{0}$ in \cite{Hansen_1.5, Hansen_2}. To avoid any
confusion we shall in the following write $\mathbf{F}_{2}^{0,new}$ for the
set denoted by $\mathbf{F}_{2}^{0}$ in \cite{Hansen_3}.

We start with a discussion of the treatment of Aitken's Theorem in \cite%
{Hansen_3}. Hansen first gives a rendition of the classical Aitken Theorem
(Theorem 3 in \cite{Hansen_3}) which is identical to Theorem 3 in \cite%
{Hansen_1.5, Hansen_2}, and thus to Theorem \ref{Aitken_variant} in the
preceding section. He proceeds to provide his `modern Aitken Theorem'
(Theorem 4 in \cite{Hansen_3}), which is identical to the corresponding
Theorem 4 in \cite{Hansen_1.5, Hansen_2} and which in turn is equivalent to
Assertion 4 as just discussed in Section \ref{Ait} above.\footnote{%
The same caveat regarding the formulation of Hansen's theorems as discussed
in Section \ref{Ait} applies here.} Consequently, the discussion given in
Section \ref{Ait} above applies. In particular, the estimators figuring in
the `modern Aitken Theorem' in \cite{Hansen_3} are all automatically linear
by our Theorem \ref{class_trivial}, and hence the `modern Aitken Theorem' in 
\cite{Hansen_3} is not new, but reduces to the classical Aitken Theorem.

Hansen then goes on to provide a further result (Theorem 5 in \cite{Hansen_3}%
) which can equivalently be stated as the following assertion (the
equivalence is seen in the same way as the equivalence between Theorem 3 in 
\cite{Hansen_1.5, Hansen_2} and Theorem \ref{Aitken_variant} in Section \ref%
{Ait} above).\footnote{%
The same caveat regarding the formulation of Hansen's theorems discussed in
Section \ref{Ait} applies also to Theorem 5 in \cite{Hansen_3}.} For $\Sigma 
$ a diagonal $n\times n$ matrix with positive diagonal elements, define $%
\mathbf{F}_{2}^{\ast }(\Sigma )=\{F\in \mathbf{F}_{2}^{\ast
}:Var_{F}(e)\propto \Sigma \}$, where $\propto $ denotes proportionality. Of
course, then $\mathbf{F}_{2}^{\ast }=\tbigcup \{\mathbf{F}_{2}^{\ast
}(\Sigma ):\Sigma $ diagonal with positive diagonal elements$\}$ and $%
\mathbf{F}_{2}^{\ast }(I_{n})=\mathbf{F}_{2}^{0,new}$. Recall that $\hat{%
\beta}_{GLS}=\hat{\beta}_{GLS}(\Sigma )=(X^{\prime }\Sigma
^{-1}X)^{-1}X^{\prime }\Sigma ^{-1}Y$.

\bigskip

\textbf{Assertion 5: }Let $\Sigma $ be an arbitrary known \emph{diagonal} $%
n\times n$ matrix with positive diagonal elements. If $\hat{\beta}$ is an
estimator that is unbiased under all $F\in \mathbf{F}_{2}^{\ast }$ (meaning
that $E_{F}\hat{\beta}=\beta (F)$ for every $F\in \mathbf{F}_{2}^{\ast }$),
then 
\begin{equation*}
Var_{F}(\hat{\beta})\succeq Var_{F}(\hat{\beta}_{GLS})
\end{equation*}%
for every $F\in \mathbf{F}_{2}^{\ast }(\Sigma )$.

\bigskip

This assertion is very much different from an Aitken Theorem:\footnote{%
Assertion 5 (equivalently, Theorem 5 in \cite{Hansen_3}) seems to be
correct. However, we have not checked the correctness of the proofs in
Section 6 of \cite{Hansen_3} in any detail.} (i) The errors in the model
corresponding to $F\in \mathbf{F}_{2}^{\ast }(\Sigma )$ (i.e., distributions 
$F$ for which the variance inequality has to hold) need to be \emph{%
independent}, an assumption alien to Aitken's Theorem (even if $\Sigma $ is
diagonal) as this theorem relies only on first and second moment assumptions
(as opposed to an independence assumption). (ii) The unbiasedness assumption
is -- like in results discussed earlier -- required to hold under the \emph{%
wider} class of distributions $\mathbf{F}_{2}^{\ast }$, and not only under
the distributions $F$ describing the data generating mechanism (i.e.,~$F\in 
\mathbf{F}_{2}^{\ast }(\Sigma )$).\footnote{\label{fn_star}Requiring
unbiasedness over the wider class $\mathbf{F}_{2}^{\ast }$ is crucial here:
If in Assertion 5 unbiasedness is only required to hold for $F\in \mathbf{F}%
_{2}^{\ast }(\Sigma )$ rather than for $F\in \mathbf{F}_{2}^{\ast }$, the
resulting statement is incorrect in general. This follows for $\Sigma =I_{n}$
from Example \ref{EX2} in Appendix \ref{App_A}. [Note that the estimator
constructed in this example is unbiased even under every $F\in \mathbf{F}%
_{2}(I_{n})=\mathbf{F}_{2}^{0}$, and that the offending distribution
constructed in this example belongs to $\mathbf{F}_{2}^{\ast }(I_{n})$, in
fact even corresponds to independent identically distributed errors with
finite second moments.] Another counterexample is provided by Example \ref%
{EX1} in Appendix \ref{App_A}, which covers the location case. Similar
examples can easily be constructed for any diagonal $\Sigma $ with positive
diagonal elements by a transformation argument.} And (iii) an Aitken Theorem
should allow for general $\Sigma $. While in the context of Assertions 2 and
4 no nonlinear unbiased estimator exists, in the context of Assertion 5
nonlinear unbiased estimators indeed exist (at least for some matrices $X$),
cf.~Remark \ref{rem_class}\ below.

We next turn to the treatment of the Gauss-Markov Theorem in \cite{Hansen_3}%
: He starts with Theorem 1 (for linear estimators), which despite given the
label Gauss-Markov, is \emph{not} the Gauss-Markov Theorem, but a much
weaker result relying on the unnecessarily restrictive assumption that the
errors in the model are \emph{independent} (and homoskedastic).\footnote{%
A similar caveat as in Footnote \ref{FN1} also applies to Theorems 1, 6, and
7 in \cite{Hansen_3}.} Such an independence assumption is superfluous in the
classical Gauss-Markov Theorem. [Note that as long as only linear estimators
are considered, requiring unbiasedness for all $F\in \mathbf{F}_{2}^{\ast }$%
, as is done in Theorem 1 of \cite{Hansen_3}, is identical to requiring
unbiasedness for all $F\in \mathbf{F}_{2}^{0,new}$, or for all $F\in \mathbf{%
F}_{2}^{0}=\mathbf{F}_{2}(I_{n})$ for that matter; cf.~the discussion
following Theorem \ref{GMHansen} in Section \ref{sect:GM}.] The `modern
Gauss-Markov Theorem' (Theorem 6 in \cite{Hansen_3}) is now the special case
of Assertion 5 for $\Sigma =I_{n}$; note that this result is just Theorem 1
in \cite{Hansen_3} with the linearity requirement dropped.\footnote{%
The `modern Gauss-Markov Theorem' as given in Theorem 5 of \ \cite%
{Hansen_1.5, Hansen_2} is no longer presented, but of course is an immediate
consequence of Theorem 4 in \cite{Hansen_3}. Also Theorem 6 of \cite%
{Hansen_1.5, Hansen_2} is no longer given.} For reasons (i) and (ii)
discussed in the preceding paragraph in connection with Assertion 5, this
result can not legitimately be called a (modern) Gauss-Markov Theorem.
Finally, Theorem 7 in \cite{Hansen_3} just specializes Theorem 6 in the same
reference to the location model, hence the same remarks apply.

To sum up, Theorems 4-7 in \cite{Hansen_3} are either an intransparent
restatement of the classical Aitken Theorem introducing linearity of the
estimators through the backdoor (Theorem 4 in \cite{Hansen_3}), or are
results modelled on the Gauss-Markov or Aitken Theorem but employing \emph{%
substantial extra} conditions such as independence assumptions,
etc.~(Theorems 5-7 in \cite{Hansen_3}). [The significance and scope of the
latter results is unclear for the reasons discussed before.] As a
consequence, the advertisements regarding dropping of the linearity
assumption made in \cite{Hansen_3} are by no means justified. In particular,
the claim made in the abstract and repeated at the end of Section 3 of \cite%
{Hansen_3}, that his theorems would show that the label "linear estimator"
can be dropped from the pedagogy of the Gauss-Markov Theorem, is without any
base. We thus repeat our warning against dropping the linearity assumption
from the Gauss-Markov or Aitken Theorem.

\begin{remark}
\label{rem_class} In the discussion following Theorem 5 in \cite{Hansen_3},
the author gives an example of a nonlinear estimator that is unbiased under
all $F\in \mathbf{F}_{2}^{\ast }$. The object$\ \tilde{\beta}$ given there,
however, is not well-defined as it is the sum of two components that are of
different dimension (unless $k=1$). This can be rectified by redefining $%
\tilde{\beta}$ as $\hat{\beta}_{OLS}+Y_{i}(Y_{j}-X_{j}^{\prime }\hat{\beta}%
_{-i})a$ ($i\neq j$) for any chosen $k\times 1$ vector $a\neq 0$ (here $%
Y_{j} $ and $X_{j}^{\prime }$ denote the $j$-th row of $Y$ and $X$,
respectively).\footnote{%
There is an implicit assumption here, namely that the design matrix
continues to have full column rank even after the $i$-th row is deleted.}
This object is indeed unbiased under all $F\in \mathbf{F}_{2}^{\ast }$ (but,
in general, \emph{not} under all $F\in \mathbf{F}_{2}$). The claim in \cite%
{Hansen_3} that this is a \emph{nonlinear} estimator, however, is not
generally true for any design matrix $X$. For example, if $n=k+1$, then any
leave-one-out residual is zero, and hence $\tilde{\beta}=\hat{\beta}_{OLS}$
is linear (for any choice of $i$ and $j$). Another example where $\tilde{%
\beta}=\hat{\beta}_{OLS}$ is when $k=1$, the regressor is the first standard
basis vector, $i\neq 1$, and $j=1$. Fortunately, there are examples of
design matrices for which $\tilde{\beta}$ is indeed truly nonlinear.
\end{remark}

\section{Independent Identically Distributed Errors\label{sect:iid}}

We round-off the discussion by briefly considering in this section what
happens if we add the condition 
\begin{equation}
e_{1},\ldots ,e_{n}\text{ \ are i.i.d.}  \label{iid}
\end{equation}%
to the model. Let $\boldsymbol{F}_{2}^{iid}$ be the subset of $\boldsymbol{F}%
_{2}^{0}$ corresponding to distributions $F$ that result from (\ref{model}),
(\ref{expect}), (\ref{var}), and (\ref{iid}). In particular, we ask what is
the status of the following assertion which is analogous to Assertion 1.

\bigskip

\textbf{Assertion 6:} If $\hat{\beta}$ is an estimator that is unbiased
under all $F\in \mathbf{F}_{2}^{iid}$ (meaning that $E_{F}\hat{\beta}=\beta
(F)$ for every $F\in \mathbf{F}_{2}^{iid}$), then 
\begin{equation*}
Var_{F}(\hat{\beta})\succeq Var_{F}(\hat{\beta}_{OLS})
\end{equation*}%
for every $F\in \mathbf{F}_{2}^{iid}$.\bigskip

Note that Assertion 6 differs from Assertion 1 in two respects: (i) the set
of competitors to $\hat{\beta}_{OLS}$, i.e., the set of unbiased estimators
in Assertion 6 is potentially larger than the corresponding set in Assertion
1, and (ii) the set of distributions $F$ for which the variance inequality
has to hold has gotten smaller compared to Assertion 1. Hence, the
truth-status of Assertion 1 does not inform us about the corresponding
status of Assertion 6.

Fortunately, Example \ref{EX2} in Appendix \ref{App_A} comes to the rescue
and shows that Assertion 6 is incorrect in general (meaning that a design
matrix can be found such that it is false). This is so since the nonlinear
estimator constructed in that example is a fortiori unbiased under $\mathbf{F%
}_{2}^{iid}$, and since the offending $F$ found in that example in fact
belongs to $\mathbf{F}_{2}^{iid}$. However, in the special case of the
location model Assertion 6 is actually true. This follows directly from
Theorem 5 in \cite{Halmos}.\footnote{\cite{Halmos} allows for $\sigma ^{2}=0$%
. However, this is immaterial as a consequence of the discussion in Footnote %
\ref{fn_koop}.} [Recall that, in contrast, Assertion 1 is false in the case
of a location model; cf. Example \ref{EX1} in Appendix \ref{App_A}.\footnote{%
It is perhaps interesting to note that the assertion one obtains from
Assertion 6 by replacing $\mathbf{F}_{2}^{iid}$ by $\mathbf{F}_{2}^{\ast
}(I_{n})$ at every occurrence in Assertion 6 is also incorrect in general,
and even in the location case; see the discussion in Footnote \ref{fn_star}.}%
]

For results in the location case pertaining to classes of absolutely
continuous distributions (without or with symmetry restrictions) see Example
4.2 in Section 2.4 of \cite{LehCas} and the discussion following this
example.

A nice result is due to \cite{KS}: Suppose we restrict to i.i.d.~errors in
our regression model, but where now the distribution of the errors, $G$ say,
is known (and has finite second moments). Suppose also that $n\geq 2k+1$ and
that the design matrix has no rows of zeroes. Then, if $\hat{\beta}_{OLS}$
is best unbiased in this model, the distribution $G$ must be Gaussian. [\cite%
{KS} actually prove a more general result.] A related result for the
location model with independent (not necessarily identically distributed)
errors is given in Theorem 7.4.1 of \cite{KLR}. For more results in that
direction see Sections 7.4-7.9 in the same reference.

There is probably more in the mathematical statistics literature we are not
aware of, but this is what a quick search has turned up.

\appendix{}

\section{Appendix: Counterexamples\label{App_A}}

Here we provide various counterexamples to Assertion 1. They all rest on the
following lemma which certainly is not original as similar computations can
be found in the literature, see, e.g., \cite{Gnotetal}, \cite{Knautz_1993,
Knautz_1999}, and references therein. Counterexamples can also be easily
derived from results in the before mentioned papers. In this appendix we
always maintain the model from Section \ref{frame} and assume that (\ref%
{homo}) holds. For the case $\Sigma \neq I_{n}$, similar counterexamples to
Assertion 3 can be obtained by a standard transformation argument. We do not
pursue this any further.

\begin{lemma}
\label{help} Consider the model as in Section \ref{frame}, additionally
satisfying (\ref{homo}).

(a) Define estimators via%
\begin{equation}
\hat{\beta}_{\alpha }=\hat{\beta}_{OLS}+\alpha (Y^{\prime }H_{1}Y,\ldots
,Y^{\prime }H_{k}Y)^{\prime }  \label{quadr_est}
\end{equation}%
where the $H_{i}$'s are symmetric $n\times n$ matrices and $\alpha $ is a
real number. Suppose $\limfunc{tr}(H_{i})=0$ and $X^{\prime }H_{i}X=0$ for $%
i=1,\ldots ,k$. Then $E_{F}(\hat{\beta}_{\alpha })=\beta (F)$ for all $F\in 
\mathbf{F}_{2}^{0}$.

(b) Suppose the $H_{i}$'s are as in Part (a). If $Cov_{F}(c^{\prime }\hat{%
\beta}_{OLS},c^{\prime }(Y^{\prime }H_{1}Y,\ldots ,Y^{\prime
}H_{k}Y)^{\prime })\neq 0$ for some $c\in \mathbb{R}^{k}$ and for some $F\in 
\mathbf{F}_{2}^{0}$ with finite fourth moments, then there exists an $\alpha
\in \mathbb{R}$ such that%
\begin{equation}
Var_{F}(c^{\prime }\hat{\beta}_{\alpha })<Var_{F}(c^{\prime }\hat{\beta}%
_{OLS});  \label{inequ}
\end{equation}%
in particular, $\hat{\beta}_{OLS}$ then does not have smallest
variance-covariance matrix (w.r.t.~Loewner order) over $\mathbf{F}_{2}^{0}$
in the class of all estimators that are unbiased under all $F\in \mathbf{F}%
_{2}^{0}$.\footnote{%
Recall the convention discussed in Remark \ref{rem_Loew}.}

(c) Suppose the $H_{i}$'s are as in Part (a). For every $c\in \mathbb{R}^{k}$
and for every $F\in \mathbf{F}_{2}^{0}$ (with finite fourth moments) under
which $\beta (F)=0$ we have%
\begin{equation}
Cov_{F}\left( c^{\prime }\hat{\beta}_{OLS},c^{\prime }(Y^{\prime
}H_{1}Y,\ldots ,Y^{\prime }H_{k}Y)^{\prime }\right)
=\tsum_{j=1}^{n}\tsum_{l=1}^{n}\tsum_{m=1}^{n}d_{j}\left(
\tsum_{i=1}^{k}c_{i}h_{lm}(i)\right) E_{F}(e_{j}e_{l}e_{m}),  \label{cov_0}
\end{equation}%
where $d=(d_{1},\ldots ,d_{n})^{\prime }=X(X^{\prime }X)^{-1}c$ and $%
h_{lm}(i)$ denotes the $(l,m)$-th element of $H_{i}$.

(d) Suppose the $H_{i}$'s are as in Part (a). For every $c\in \mathbb{R}^{k}$
and for every $F\in \mathbf{F}_{2}^{0}$ (with finite fourth moments) under
which (i) $\beta (F)=0$ and under which (ii) the coordinates of $Y$ are
independent (equivalently, the errors $e_{i}$ are independent)%
\begin{equation}
Cov_{F}\left( c^{\prime }\hat{\beta}_{OLS},c^{\prime }(Y^{\prime
}H_{1}Y,\ldots ,Y^{\prime }H_{k}Y)^{\prime }\right)
=\tsum_{j=1}^{n}d_{j}\left( \tsum_{i=1}^{k}c_{i}h_{jj}(i)\right)
E_{F}(e_{j}^{3}).  \label{cov}
\end{equation}
\end{lemma}

\textbf{Proof:} The proof of Parts (a), (c), and (d) is by straightforward
computation. Since 
\begin{eqnarray}
Var_{F}(c^{\prime }\hat{\beta}_{\alpha }) &=&Var_{F}(c^{\prime }\hat{\beta}%
_{OLS})+2\alpha Cov_{F}\left( c^{\prime }\hat{\beta}_{OLS},c^{\prime
}(Y^{\prime }H_{1}Y,\ldots ,Y^{\prime }H_{k}Y)^{\prime }\right)  \notag \\
&&+\alpha ^{2}Var_{F}(c^{\prime }(Y^{\prime }H_{1}Y,\ldots ,Y^{\prime
}H_{k}Y)^{\prime }),  \label{Var}
\end{eqnarray}%
the claim in (b) follows immediately as the first derivative of $%
Var_{F}(c^{\prime }\hat{\beta}_{\alpha })$ w.r.t. $\alpha $ and evaluated at 
$\alpha =0$ equals $2Cov_{F}\left( c^{\prime }\hat{\beta}_{OLS},c^{\prime
}(Y^{\prime }H_{1}Y,\ldots ,Y^{\prime }H_{k}Y)^{\prime }\right) $. Note that
all terms in (\ref{Var}) are well-defined and finite because of our fourth
moment assumption. Hence, whenever this covariance is non-zero, we may
choose $\alpha \neq 0$ small enough such that (\ref{inequ}) holds. $%
\blacksquare $

We now provide a few counterexamples that make use of the preceding lemma.

\begin{example}
\label{EX1}Consider the location model, i.e., the case where $k=1$ and $%
X=(1,\ldots ,1)^{\prime }$. Choose $H_{1}$ as the $n\times n$ matrix which
has $h_{11}(1)=-h_{22}(1)=1$ and $h_{ij}(1)=0$ else. Then the conditions on $%
H_{1}$ in Part (a) of Lemma \ref{help} are satisfied, and hence $\hat{\beta}%
_{\alpha }$ is unbiased under all $F\in \mathbf{F}_{2}^{0}$. Setting $c=1$,
we find for the covariance in (\ref{cov})%
\begin{equation*}
n^{-1}(E_{F}(e_{1}^{3})-E_{F}(e_{2}^{3}))\neq 0
\end{equation*}%
for every $F$ $\in \mathbf{F}_{2}^{0}$ (with finite fourth moments) under
which $\beta (F)=0$, the errors $e_{i}$ are independent, and $%
E_{F}(e_{1}^{3})\neq E_{F}(e_{2}^{3})$ hold. Such distributions $F$
obviously exist.\footnote{%
E.g., choose $e_{2},\ldots ,e_{n}$ i.i.d.$~N(0,\sigma ^{2})$ and $e_{1}$
independent from $e_{2},\ldots ,e_{n}$ with mean zero, variance $\sigma ^{2}$%
, third moment not equal to zero, and finite fourth moment.} As a
consequence, $\hat{\beta}_{OLS}$ is not best (over $\mathbf{F}_{2}^{0}$) in
the class of all estimators $\hat{\beta}$ that are unbiased under all $F\in 
\mathbf{F}_{2}^{0}$. In particular, Assertion 1 is false for this design
matrix.
\end{example}

For the argument underlying the preceding example it is key that the errors
are \emph{not} i.i.d.~under the relevant $F$. In fact, in the location model
(i.e., $X=(1,\ldots ,1)^{\prime }$) we have $Var_{F}(\hat{\beta}%
_{OLS})\preceq Var_{F}(\hat{\beta}_{\alpha })$ for every real $\alpha $, for
every choice of $H_{1}$ as in Part (a) of Lemma \ref{help}, and for every $%
F\in \mathbf{F}_{2}^{0}$ (with finite fourth moments) under which the errors 
$e_{i}$ are i.i.d., since then $Cov_{F}(\hat{\beta}_{OLS},Y^{\prime
}H_{1}Y)=0$ as is easily seen. [This is in line with the result of \cite%
{Halmos} discussed in Section \ref{sect:iid}.] For other design matrices $X$
the argument, however, works even for i.i.d.~errors as we show in the
subsequent example. Cf.~Section 4.1 of \cite{Gnotetal} for related results
and more.

\begin{example}
\label{EX2}Consider the balanced one-way layout for $k=2$ and $n=4$. That
is, $X$ has first column equal to $(1,1,0,0)^{\prime }$ and second column
equal to $(0,0,1,1)^{\prime }$. Set $c=(1,0)^{\prime }$. Then $%
d=(1/2,1/2,0,0)^{\prime }$. Choose, e.g., $H_{1}=H_{2}$ as the $4\times 4$
matrix made up of $2\times 2$ blocks, where the off-diagonal blocks are
zero, the first and second diagonal block, respectively, are given by 
\begin{equation*}
\left( 
\begin{array}{cc}
1 & -1 \\ 
-1 & 1%
\end{array}%
\right) ,\left( 
\begin{array}{cc}
-1 & 1 \\ 
1 & -1%
\end{array}%
\right) .
\end{equation*}%
Then the conditions on $H_{i}$ in Part (a) of Lemma \ref{help} are
satisfied, and hence $\hat{\beta}_{\alpha }$ is unbiased under all $F\in 
\mathbf{F}_{2}^{0}$. For the covariance in (\ref{cov}) we find%
\begin{equation*}
(E_{F}(e_{1}^{3})+E_{F}(e_{2}^{3}))/2
\end{equation*}%
under any $F$ $\in \mathbf{F}_{2}^{0}$ (with finite fourth moments) under
which $\beta (F)=0$ and the errors $e_{i}$ are independent. If $F$ is chosen
such that the errors are furthermore i.i.d.~and asymmetrically distributed,
the expression in the preceding display reduces to $E_{F}(e_{1}^{3})\neq 0$.
Such distributions $F$ obviously exist. As a consequence, $\hat{\beta}_{OLS}$
is not best (over $\mathbf{F}_{2}^{0}$) in the class of all estimators $\hat{%
\beta}$ that are unbiased under all $F\in \mathbf{F}_{2}^{0}$. In
particular, Assertion 1 is false for this design matrix.
\end{example}

Many more counterexamples can be generated with the help of Lemma \ref{help}
as outlined in the subsequent remark.

\begin{remark}
(i) Suppose $X$ admits a choice of $H_{i}$ satisfying the conditions in Part
(a) of Lemma \ref{help} and a $c\in \mathbb{R}^{k}$ such that $%
\tsum_{j=1}^{n}d_{j}\tsum_{i=1}^{k}c_{i}h_{jj}(i)\neq 0$. Then the
covariance in (\ref{cov}) is not zero if $F$ in Part (d) of the lemma is
chosen to correspond to asymmetrically distributed i.i.d.~errors. Part (b)
of the lemma can then be applied. In case $H_{i}=H$ for all $i=1,\ldots ,k$,
these conditions further reduce to $\tsum_{j=1}^{n}d_{j}h_{jj}\neq 0$ and $%
\tsum_{i=1}^{k}c_{i}\neq 0$.

(ii) Suppose $X$ admits a choice of $H_{i}$ satisfying the conditions in
Part (a) of Lemma \ref{help} and a $c\in \mathbb{R}^{k}$ such that for an
index $j_{0}$ it holds that $d_{j_{0}}\tsum_{i=1}^{k}c_{i}h_{j_{0}j_{0}}(i)%
\neq 0$. Then the covariance in (\ref{cov}) is not zero if $F$ in Part (d)
of the lemma is chosen to correspond to independent errors with $%
E_{F}(e_{j_{0}}^{3})\neq 0$ and $E_{F}(e_{j}^{3})=0$ for $j\neq j_{0}$.
Again Part (b) of the lemma can then be applied. In case $H_{i}=H$ for all $%
i=1,\ldots ,k$, these conditions further reduce to $d_{j_{0}}h_{j_{0}j_{0}}%
\neq 0$ and $\tsum_{i=1}^{k}c_{i}\neq 0$.

(iii) Part (c) of Lemma \ref{help} allows for further examples to be
generated, where now the errors need not be independently distributed under
the relevant $F$.
\end{remark}

One certainly could set out to characterize those design matrices $X$ for
which a counterexample to Assertion 1 can be constructed with the help of
Lemma \ref{help}. We do not pursue this here. In particular, we have not
investigated whether for \emph{any} $n\times k$ design matrix $X$ with $k<n$
one can construct an estimator $\hat{\beta}_{\alpha }$ as in the lemma that
satisfies (\ref{inequ}) for some $c\in \mathbb{R}^{k}$ and for some $F\in 
\mathbf{F}_{2}^{0}$.

\section{Appendix: Adding A Finite Second Moment Assumption on the
Estimators \label{App_B}}

We start by discussing the consequences of introducing the requirement $%
E_{F}(\parallel \hat{\beta}\parallel ^{2})<\infty $ for every $F\in \mathbf{F%
}_{2}^{0}$ into Assertions 1 and 2. First, note that nothing changes for
Assertion 2 (and for the corresponding Theorem 5 in \cite{Hansen_1.5,
Hansen_2}), since the estimators $\hat{\beta}$ allowed in Assertion 2 are
all linear as a consequence of our Theorem \ref{class_trivial} (and thus
have finite second moments even for all $F\in \mathbf{F}_{2}$). Turning to
Assertion 1, we first give the following proposition.

\begin{proposition}
Consider estimators of the form%
\begin{equation*}
\hat{\beta}=AY+(Y^{\prime }H_{1}Y,\ldots ,Y^{\prime }H_{k}Y)^{\prime }
\end{equation*}%
with $AX=I_{k}$ and the $n\times n$ matrices $H_{j}$ satisfying $\limfunc{tr}%
(H_{j})=0$ and $X^{\prime }H_{j}X=0$ for $j=1,\ldots ,k$. Suppose $%
E_{F}(\parallel \hat{\beta}\parallel ^{2})<\infty $ for every $F\in \mathbf{F%
}_{2}^{0}$ holds. Then $\hat{\beta}$ is a linear estimator.
\end{proposition}

\textbf{Proof:} First, observe that $E_{F}(\hat{\beta})=\beta (F)$ for all $%
F\in \mathbf{F}_{2}^{0}$ holds. Second, we may assume the matrices $H_{j}$
to be symmetric (if necessary we replace $H_{j}$ by $(H_{j}+H_{j}^{\prime
})/2$). Since $AY$ obviously has finite second moments under every $F\in 
\mathbf{F}_{2}^{0}$, the finite second moment assumption on $\hat{\beta}$
implies that $E_{F}((Y^{\prime }H_{j}Y)^{2})<\infty $ has to hold for every $%
F\in \mathbf{F}_{2}^{0}$ and every $j=1,\ldots ,k$. Let now $j$ ($j=1,\ldots
,k$) be arbitrary, but fixed. By symmetry of $H_{j}$, there exists an
orthogonal matrix $U$ such that $UH_{j}U^{\prime }=\Lambda _{j}$ where $%
\Lambda _{j}$ is a diagonal matrix. Set $Z=UY$ and let $\lambda _{ij}$
denote the elements on the diagonal of $\Lambda _{j}$. Then%
\begin{equation}
(Y^{\prime }H_{j}Y)^{2}=(Z^{\prime }\Lambda _{j}Z)^{2}=\left(
\tsum_{i=1}^{n}\lambda _{ij}z_{i}^{2}\right) ^{2}=\tsum_{i=1}^{n}\lambda
_{ij}^{2}z_{i}^{4}+\tsum_{i,l=1,i\neq l}^{n}\lambda _{ij}\lambda
_{lj}z_{i}^{2}z_{l}^{2}.  \label{quartic}
\end{equation}%
Let $F\in \mathbf{F}_{2}^{0}$ be such that $\beta (F)=0$ and such that the
elements of $Z$ are i.i.d.~with mean zero, finite variance, and infinite
fourth moment. Such an $F$ exists: Start with a distribution on $Z$ with the
required properties just listed and work backwards, defining $Y=U^{\prime }Z$%
. Then clearly the implied $F$ has $\beta (F)=0$ and belongs to $\mathbf{F}%
_{2}^{0}$. Since the coordinates of $Z$ are independent and have a finite
second moment, the term $\tsum_{i,l=1,i\neq l}^{n}\lambda _{ij}\lambda
_{lj}z_{i}^{2}z_{l}^{2}$ has finite expectation under $F$. Since the
l.h.s.~of (\ref{quartic}) has finite expectation under $F$ under our
assumptions as noted before, $E_{F}(\tsum_{i=1}^{n}\lambda
_{ij}^{2}z_{i}^{4})$ has to be finite. Since $E_{F}(z_{i}^{4})=\infty $ for
every $i=1,\ldots ,n$, we must have $\lambda _{ij}=0$ for every $i$ and the
given $j$. This shows that $H_{j}=0$. Since $j$ was arbitrary, this holds
for every $j$, and thus $\hat{\beta}=AY$ is linear. $\blacksquare $

The significance of the preceding proposition is the following. Suppose
Theorem 4.3 in \cite{koopmann} is indeed correct (recall that no complete
proof is given in \cite{koopmann}), and thus only quadratic estimators as in
the above proposition figure in Assertion 1.\footnote{%
Cf.~also Footnote \ref{fn_koop}.} Then introducing the \emph{extra}
condition of finite second moments for the estimators $\hat{\beta}$ (under 
\emph{every} $F\in \mathbf{F}_{2}^{0}$) into Assertion 1 reduces the class
of competitors to the class of linear unbiased estimators. The resulting 
\emph{version} of Assertion 1 is thus true and coincides with the classical
Gauss-Markov Theorem. [Recall that Assertion 1 is false in general.] Again,
linearity is reintroduced indirectly by adding the before mentioned extra
condition. In case Theorem 4.3 in \cite{koopmann} is incorrect, then other
estimators than quadratic ones might figure in this version of Assertion 1
and it is unclear whether this version of Assertion 1 is true or not, and if
it is true what its scope is.

Concerning Assertion 3, an analogous version of the preceding proposition
can be given and a similar discussion applies. Note that Assertion 4 and the
corresponding theorems in Hansen's papers are not affected by introducing an
extra finite second moment assumption on the estimators, because all
estimators in Assertion 4 are automatically linear in view of our Theorem %
\ref{class_trivial} (and thus have finite second moments under all $F\in 
\mathbf{F}_{2}$).

Adding an extra finite second moment assumption on the estimators to
Assertion 5 leads to a version that, a fortiori, is not an `Aitken Theorem',
for the same reasons as discussed in Section \ref{H3}. [In case Assertion 5
is correct, this version is a fortiori also correct.] Whether or not
Assertion 6, which is false in general, becomes a true statement after
adding an extra finite second moment condition on the estimators, we have
not investigated.

\bibliographystyle{ims}
\bibliography{refs}

\end{document}